\documentclass{article}
\usepackage{graphicx} 
\usepackage[utf8]{inputenc}
\usepackage{amsmath,amssymb,amsthm, amsfonts}
\usepackage[margin=1in]{geometry}
\usepackage{enumitem}
\usepackage[colorlinks=true, linkcolor=blue, citecolor=blue, urlcolor=blue]{hyperref}
\usepackage[capitalise, nameinlink]{cleveref}
\usepackage{tikz}
\usepackage{float}

\newtheorem{theorem}{Theorem}[subsection]
\newtheorem{lemma}[theorem]{Lemma}
\newtheorem{proposition}[theorem]{Proposition}
\newtheorem{definition}[theorem]{Definition}
\newtheorem{example}[theorem]{Example}

\newcommand{\N}{\mathbb{N}}
\newcommand{\rev}{\operatorname{rev}}
\newcommand{\hgt}{\operatorname{ht}}
\newcommand{\lowid}[1]{#1\!\downarrow}
\newlist{conditions}{enumerate}{1}
\setlist[conditions]{label=(C\arabic*), ref=C\arabic*}
\crefname{conditionsi}{condition}{conditions}

\newcommand{\pge}[1]{\mathrel{\mathchoice{\ge_{#1}}{\ge_{#1}}{\,\ge_{#1}\,}{\,\ge_{#1}\,}}}
\newcommand{\pgt}[1]{\mathrel{\mathchoice{>_{#1}}{>_{#1}}{\,>_{#1}\,}{\,>_{#1}\,}}}

\newcommand{\va}{\mathbf{a}}
\newcommand{\vb}{\mathbf{b}}
\newcommand{\vc}{\mathbf{c}}
\newcommand{\vk}{\mathbf{k}}
\newcommand{\vm}{\mathbf{m}}

\newcommand{\vx}{\mathbf{x}}
\newcommand{\vy}{\mathbf{y}}

\title{The reversal symmetry of height in the stable Tamari lattice}
\author{Kiran Mathew}
\date{}
\begin{document}

\maketitle

\begin{abstract}
The stable Tamari lattice, introduced by Chau et al., is a partial ordering of sequences of nonnegative integers. We study a particular symmetry of the height generating function $f_{\vb}(t)$, which gives the height distribution for all intervals in the lattice with a fixed top element $\vb$. It has been conjectured that the function is invariant under reversal of this top sequence. We associate each element under $\vb$ with a unique ID vector, whose weight is related to interval height by a simple formula. We then construct a height-preserving bijection from ID vectors under $\vb$ to ID vectors under the reversal of $\vb$ to show the symmetry of the height generating functions.
\end{abstract}

\section{Introduction}\label{sec-1}
Bergeron and Pr\'eville-Ratelle's study of diagonal harmonics led them to introduce the $r$-Tamari lattice, an extension of the classical Tamari lattice~\cite{MR3029440}. Within it, they studied the interval generating functions
\[
  i^{(r)}_\beta(q) = \sum_{\alpha \le \beta} q^{d(\alpha,\beta)},
\]
where $d(\alpha, \beta)$ is the length of a longest chain from $\alpha$ to $\beta$, as part of their conjectures on higher trivariate diagonal harmonics. Further results on the enumeration of $r$-Tamari and generalized Tamari intervals appear in~\cite{MR2880681, MR3632566}.

Chau, Daugherty, Grate, Haiman, Mart\'inez Mori, and Pun later introduced the \emph{stable Tamari lattice}~\cite{chau2026stabletamarilattice}, an infinite lattice on sequences of $n$ nonnegative integers in which the $r$-Tamari lattices appear as a nested family of principal lower ideals. A central feature of this lattice is an explicit formula for interval heights, which has encouraged study of the \emph{height generating function}
\[
  f_{\vb}(t) = \sum_{\va \le_T \vb} t^{\hgt(\va,\vb)},
\]
which records the heights of all intervals with top element $\vb$. When $\vb$ corresponds to an $r$-Dyck path $\beta$, one has $f_{\vb}(t) = i^{(r)}_\beta(t)$, so the generating functions $f_{\vb}(t)$ are stable analogues of the Bergeron--Pr\'eville-Ratelle interval generating functions.

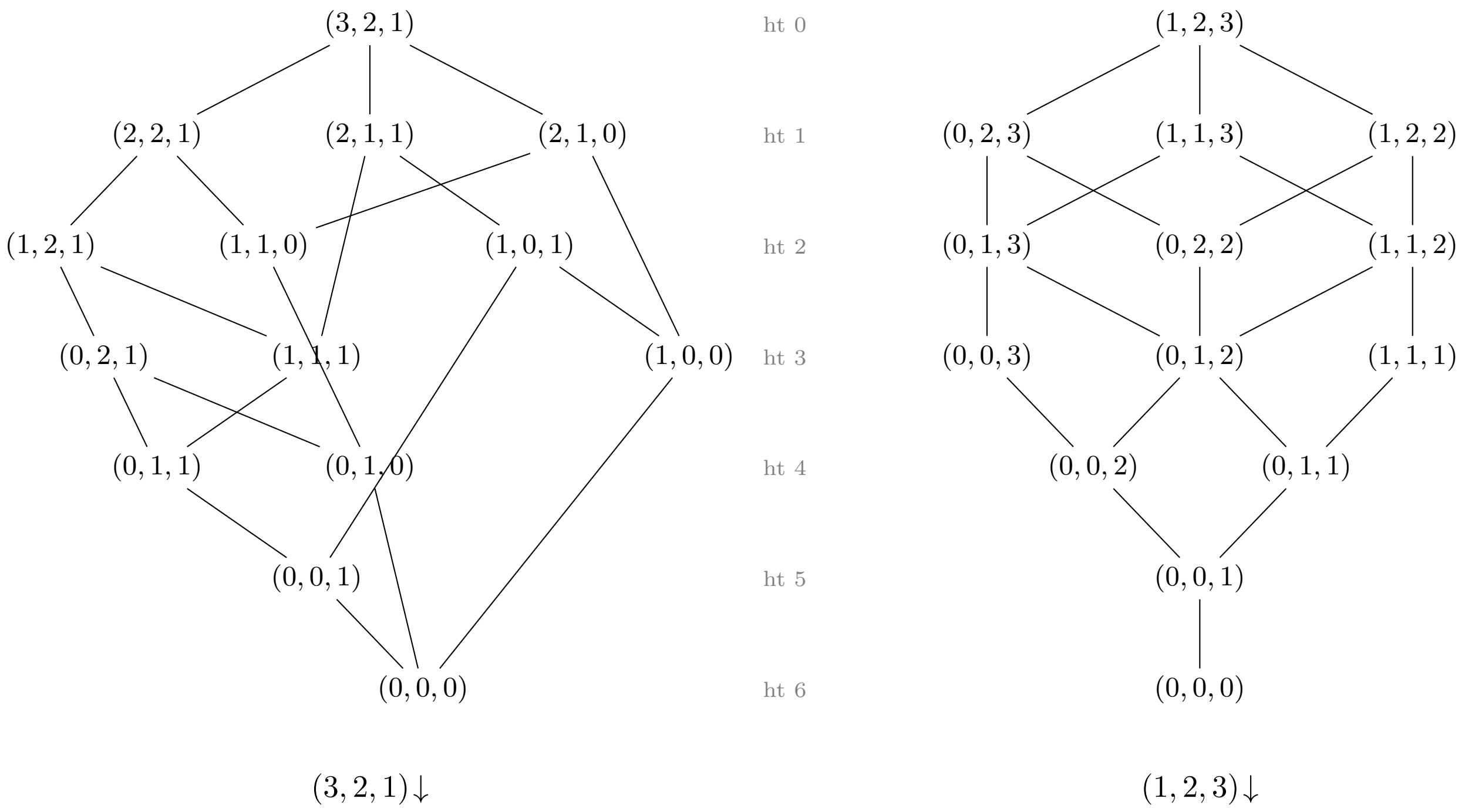
\begin{figure}[H]
\centering
\resizebox{\textwidth}{!}{%
\begin{tikzpicture}[x=1.2cm, y=1.25cm,
  every node/.style={font=\small, inner sep=2pt}]

\foreach \h in {0,...,6}
  \node[gray, font=\scriptsize] at (6.9, 6-\h) {ht $\h$};

\begin{scope}
\node (a321) at (3,   6) {$(3,2,1)$};
\node (a221) at (1,   5) {$(2,2,1)$};
\node (a211) at (3,   5) {$(2,1,1)$};
\node (a210) at (5,   5) {$(2,1,0)$};
\node (a121) at (0,   4) {$(1,2,1)$};
\node (a110) at (2,   4) {$(1,1,0)$};
\node (a101) at (4.5, 4) {$(1,0,1)$};
\node (a021) at (0.5, 3) {$(0,2,1)$};
\node (a111) at (2.5, 3) {$(1,1,1)$};
\node (a100) at (6,   3) {$(1,0,0)$};
\node (a011) at (1,   2) {$(0,1,1)$};
\node (a010) at (3,   2) {$(0,1,0)$};
\node (a001) at (2.5, 1) {$(0,0,1)$};
\node (a000) at (3.5, 0) {$(0,0,0)$};
\foreach \x/\y in {
  a210/a321, a211/a321, a221/a321,
  a100/a210, a110/a210,
  a101/a211, a111/a211,
  a110/a221, a121/a221,
  a001/a101, a100/a101,
  a010/a110,
  a021/a121, a111/a121,
  a010/a021, a011/a021,
  a011/a111,
  a000/a100, a000/a010, a001/a011, a000/a001}
  \draw (\x) -- (\y);
\node[font=\normalsize] at (3, -0.9) {$\lowid{(3,2,1)}$};
\end{scope}

\begin{scope}[shift={(7.8,0)}]
\node (c123) at (3, 6) {$(1,2,3)$};
\node (c023) at (1, 5) {$(0,2,3)$};
\node (c113) at (3, 5) {$(1,1,3)$};
\node (c122) at (5, 5) {$(1,2,2)$};
\node (c013) at (1, 4) {$(0,1,3)$};
\node (c022) at (3, 4) {$(0,2,2)$};
\node (c112) at (5, 4) {$(1,1,2)$};
\node (c003) at (1, 3) {$(0,0,3)$};
\node (c012) at (3, 3) {$(0,1,2)$};
\node (c111) at (5, 3) {$(1,1,1)$};
\node (c002) at (2, 2) {$(0,0,2)$};
\node (c011) at (4, 2) {$(0,1,1)$};
\node (c001) at (3, 1) {$(0,0,1)$};
\node (c000) at (3, 0) {$(0,0,0)$};
\foreach \x/\y in {
  c023/c123, c113/c123, c122/c123,
  c013/c023, c022/c023,
  c013/c113, c112/c113,
  c022/c122, c112/c122,
  c003/c013, c012/c013,
  c012/c022,
  c012/c112, c111/c112,
  c002/c003,
  c002/c012, c011/c012,
  c011/c111,
  c001/c002, c001/c011,
  c000/c001}
  \draw (\x) -- (\y);
\node[font=\normalsize] at (3, -0.9) {$\lowid{(1,2,3)}$};
\end{scope}
\end{tikzpicture}%
}
\caption{The posets $(3, 2, 1)\downarrow$ and $(1, 2, 3) \downarrow$ respectively, organized by height.}
\label[figure]{fig:321-vs-123}
\end{figure}

Computations led Chau et al.\ to conjecture that these generating functions are invariant under the reversal of the entries of $\vb$. For example,
\[
  f_{(3,2,1)}(t) = f_{(1,2,3)}(t) = 1 + 3t + 3t^2 + 3t^3 + 2t^4 + t^5 + t^6,
\]
even though the lower ideals $\lowid{(3,2,1)}$ and $\lowid{(1,2,3)}$ are not isomorphic as posets. Their Hasse diagrams are given in \cref{fig:321-vs-123} to illustrate this. Although the two posets have the same number of elements at each height, their covering relations are different.

Chau et al.\ construct a bijection to show that the lower ideals of $\vb$ and its reverse, $\rev(\vb)$, are of the same size \cite[Theorem~3.20]{chau2026stabletamarilattice}; equivalently, the height generating functions agree when $t=1$. However, the bijection does not preserve height between the lower ideals, meaning the full reversal symmetry \cite[Conjecture~3.12]{chau2026stabletamarilattice} remained conjectural.

In this paper, we provide the symmetry in \cref{main-result}. We define ID vectors, unique vector representatives of lower-ideal elements, for which interval height is given by $|\mathbf{b}|-|\mathbf{k}|$. We then characterize them and provide a weight-preserving bijection to connect the lower ideals of $\vb$ and $\rev(\vb)$ while preserving height. Along the way, the characterization of ID vectors gives an explicit, reconstruction-free description of any given lower ideal.

The paper is organized as follows. \cref{sec-2} recalls the necessary definitions and results on the stable Tamari lattice from Chau et al.\ and proves the characterization of ID vectors. \cref{sec-3} introduces a convenient framework for the bijection, then establishes the bijective map and deduces from it the reversal symmetry of $f_{\vb}(t)$.

\section{The Stable Tamari Order}\label{sec-2}
\subsection{Background}
We begin by recalling the definitions and results from~\cite{chau2026stabletamarilattice} that will be used throughout the paper. Boldface lowercase letters like $\va$ denote sequences of $n$ nonnegative integers, or $\N^n$. Throughout, ``sequence'' is used interchangeably with ``vector''. For a sequence $\va$, the \emph{position order} $([n], \leq_{\va})$ is defined by 
\[
i \leq_{\va} j \Longleftrightarrow i \le j \text{  and  } a_i, a_{i+1}, \ldots a_{j-1} > a_j.
\]
This is a partial order on $[n]$, the positions within a given sequence $\va$. We note that $i\le_\va i$ always, and we write $r >_{\va} i$ when specifying that $r \neq i$. We also define the \emph{weight} of $\va$ as $|\va| = \sum_{i} a_i$.

Next, we introduce the \emph{stable Tamari ordering} $(\mathbb{N}^n, \leq_T)$, which is defined by
$$\va \leq_T \vb \Longleftrightarrow \va \leq \vb \quad \text{(coordinate-wise) and} \quad a_i-a_j \leq b_i-b_j \; \forall \; i \leq_{\va} j.$$
This relation is a partial order, and $(\N^n, \le_T)$ is a lattice, namely the stable Tamari lattice. The relation is illustrated below.

\begin{example}
We consider the stable Tamari relations $\va \le_T \vb$ and $\va \le_T \vc$ for $\; \va = (3,5,1,2,1), \; \vb = (7,7,1,4,3), \; \vc = (5,5,3,3,1).$
Both $\va \le \vb$ and $\va \le \vc$ hold coordinate-wise, since $a_i \le b_i$ and $a_i \le c_i$ for all $i \in [1, 5]$. Thus, in each case only
the difference condition $a_i - a_j \le b_i - b_j$ (resp. $\le c_i - c_j$)
remains, and only over the pairs $i \le_{\va} j$ with $i < j$. These are $\;1 \le_{\va} 3, \; 2 \le_{\va} 3,$ and $\; 4 \le_{\va} 5,$
the only $i < j$ with $a_i, \ldots, a_{j-1} > a_j$; as an example, $1 \le_{\va} 3$ because $a_1,a_2>a_3$ or $3, 5 > 1$. For all of these pairs, $a_i-a_j \le b_i-b_j$, so $\va \le_T \vb$. For $\vc$, however, $a_2 - a_3 = 4 > 2 = c_2 - c_3$, contradicting the necessary condition. Hence, $\va \nleq_T \vc$, even though $\va \le \vc$ coordinate-wise.
\end{example}

Our focus is the principal lower ideal of any given sequence $\vb$, defined as $\lowid{\vb} \;= \{\va: \va \le_T \vb\}$. More specifically, we are concerned with the height of intervals of the form $[\va, \vb]$ within these lower ideals: for some $\va \in \vb \downarrow$, the height of the interval $[\va, \vb]$, or $\hgt(\va, \vb)$, is defined as the length of the longest chain from $\va$ to $\vb$. Chau et al.\ compute these heights using the matrix $M_{\va}$, which we call the \emph{height matrix}, defined for a fixed $\va$ by

\[
M_{\boldsymbol{a}} = (M_{ij}), \qquad
M_{ij} =
\begin{cases}
1, & \text{if } j \leq_{\va} i,\\
0, & \text{otherwise}
\end{cases}.
\]
Since $j \le_{\va} i$ forces $j \le i$, the matrix $M_{\va}$ is unit lower triangular, hence invertible. For $\va, \vb \in \N^n$, set
\begin{equation}
\label{m}
\vm_{\va, \vb} = (\vb-\va)M_{\va}^{-1}.
\end{equation}
Chau et al.\ show that $\va \le_T \vb$ if and only if $\vm_{\va,\vb} \in \N^n$, and that in this case,
\begin{equation}
\label{m_height}
|\vm_{\va, \vb}| = \hgt(\va, \vb).
\end{equation}
Moreover, for a fixed $\vb$, the map $\va \mapsto \vm_{\va,\vb}$ is injective on $\lowid{\vb}$, so $\va$ is determined by $\vm_{\va,\vb}$. This vector will be crucial in the later construction of the height-preserving bijection, as it gives a unique and convenient encoding of height for each $\va$.

\subsection{Combinatorial characterization of ID vectors}

We now introduce and later classify the vector that will play an important role in the bijection. For $\va \in \lowid{\vb}$, define its \emph{ID vector} by $\vk_{\va,\vb}=\vb-\vm_{\va,\vb}$. This will be referred to as an ID vector under $\vb$. We utilize these vectors to deduce information about $\hgt(\va, \vb)$, since \eqref{m_height} gives

\begin{equation}
\label{k_height}
\hgt(\va, \vb)= |\vm_{\va, \vb}| = |\vb-\vk_{\va, \vb}|=|\vb|-|\vk_{\va, \vb}|.
\end{equation}
ID vectors thereby make the relationship between a given element $\va$ and its associated $\hgt(\va, \vb)$ more apparent, and easier to manipulate. Specifically, in a bijection from an ID vector under $\vb$ to some ID vector under $\rev(\vb)$, since the weight of $\vb$ is the same as that of $\rev(\vb)$, we need only preserve the weight of the ID vectors to preserve height. We use $\vk_{\va, \vb}$ rather than $\vm_{\va, \vb}$ because the characterization is naturally expressed in terms of $\vb -\vm_{\va, \vb}$.

We first provide the direct relation between ID vectors and their associated $\va$ and $\vb$, allowing their direct construction without the intermediate step of $\vm_{\va, \vb}$ or the height matrix. Since row $r$ of the height matrix is the indicator vector of the down-set $\{\, i : i \le_{\va} r \,\}$, \eqref{m} gives
\[
  b_i - a_i = \sum_{r \pge{\va} i} m_r .
\]
Therefore, with $k_i = b_i - m_i$,
\begin{equation}
\label{k_breakdown}
  a_i = k_i - \sum_{r \pgt{\va} i} m_r = k_i + \sum_{r \pgt{\va} i} (k_r - b_r).
\end{equation}
Recall that $\le_{\va}$ is a partial order, so in particular it is transitive: $i \le_{\va} r$ and $r \le_{\va} s$ imply $i \le_{\va} s$. We also use that it is \emph{left-closed}: if $q \le_{\va} r$ and $q \le j \le r$, then $j \le_{\va} r$. We use these facts in the lemmas below. However, we first demonstrate the ID vector through the following example.

\begin{example}
Throughout we illustrate the constructions with
\[
  \vb = (3,5,2,4,1,4,3), \qquad \va = (0,4,1,3,1,4,3) \in \lowid{\vb}.
\]
It follows from \eqref{m} that $\vm_{\va, \vb} = (3,0,1,1,0,0,0)$, so $\hgt(\va,\vb) = \sum_i m_i = 5$, and
the ID vector is $\vk_{\va, \vb} = \vb - \vm_{\va, \vb} = (0,5,1,3,1,4,3).$
Beyond the reflexive relations $i \le_{\va} i$ the only order relations in $\va$ are $2 \le_{\va} 3, \quad 4 \le_{\va} 5,$ and $6 \le_{\va} 7$.
\eqref{k_breakdown} reads, for example, $k_2 = a_2 + m_3 = 4 + 1 = 5$, since $r = 3$ is the only index with $r \pgt{\va} 2$.
\end{example}

We next show that ID vectors under $\vb$ are precisely the vectors $\vx$ satisfying the following two conditions:

\begin{conditions}
  \item\label{cond:C1} $0 \le \vx \le \vb$;
  \item\label{cond:C2} whenever $i < j$ and $b_r \geq b_j$ for all $i < r < j$, either $x_i \le x_j \text{ or } x_i > b_j.$
\end{conditions}
The following two lemmas will be useful in showing that constructing the ID vector from $\va$, through \eqref{k_breakdown}, results in the expression of the above conditions. Specifically, the lemmas extract characteristics of ID vectors from the position order of the associated $\va$.

\begin{lemma}\label[lemma]{lem:k-greater-b}
If $i < j$ and $i \le_{\va} j$, then $k_i > b_j$.
\end{lemma}

\begin{proof}
Since $i \le_{\va} j$, we have $a_i > a_j$. If $r \ge_{\va} j$, then transitivity gives $r \pgt{\va} i$. Hence \eqref{k_breakdown} gives
\[
  k_i \ge a_i + \sum_{r \pge{\va} j} m_r = a_i - a_j + b_j > b_j .
\]
\end{proof}

\begin{lemma}\label[lemma]{lem:left-boundary}
Let $j \in [n]$, and suppose some index $q < j$ satisfies $q \nleq_{\va} j$.
If $q$ is the largest such index, then $k_q \le k_j$.
\end{lemma}

\begin{proof}
By the maximality of $q$, every $s$ with $q < s < j$ satisfies $s \le_{\va} j$, so
$a_{q+1}, \ldots, a_{j-1} > a_j$. Since $q \nleq_{\va} j$, this forces $a_q \le a_j$.
Next, no $r$ with $q < r \le j$ satisfies $q \le_{\va} r$: for $q < r < j$ we have
$r \le_{\va} j$, so $q \le_{\va} r$ would give $q \le_{\va} j$ by transitivity, and
$r = j$ is excluded directly. Hence every $r \pgt{\va} q$ satisfies $r > j$, and then
$r \pgt{\va} j$ by left-closure, since $q < j < r$. As $\vm_{\va,\vb} \ge 0$, this gives
\[
  \sum_{r \pgt{\va} q} m_r \;\le\; \sum_{r \pgt{\va} j} m_r ,
\]
and together with $a_q \le a_j$, \eqref{k_breakdown} yields $k_q \le k_j$.
\end{proof}

The following two lemmas will help to show that a candidate vector $\vx$ with conditions \ref{cond:C1} and \ref{cond:C2} is a valid ID vector for a  unique $\va \in \lowid{\vb}$. The first lemma shows that the position order within such $\vx$ provides additional information about the sequence. The second lemma applies this result to obtain a telescoping bound, which is later used to obtain information about recurring summation terms.

\begin{lemma}\label[lemma]{claim}
Suppose $\vx$ satisfies conditions \ref{cond:C1} and \ref{cond:C2}. Then, if $i \le_{\vx} j$, then $x_s > b_j$ for every $s$ with $i \le s < j$.
\end{lemma}

\begin{proof}
Since $i \le_{\vx} j$, we have $x_s, x_{s+1}, \ldots, x_{j-1} > x_j$ for every $s$
with $i \le s < j$; in particular $x_s > x_j$. We prove $x_s > b_j$ by downward
induction on $s$, from $s = j-1$ down to $s = i$.

In the base case of $s = j - 1$, no intermediate $r$ with $j-1 < r < j$ exist so \ref{cond:C2} forces either $x_{j-1} \le x_j$ or $x_{j-1} > b_j$, and $x_{j-1} > x_j$ eliminates the first alternative and forces $x_{j-1} > b_j$. Now fix $s$ and consider the pair $(s, j)$. Every intermediate index $r$ with
$s < r < j$ satisfies $b_r \ge x_r > b_j$ by the induction hypothesis and the condition $\vb \ge \vx$. Hence $b_r \ge b_j$ for all $s < r < j$, so
condition~\ref{cond:C2} applies to $(s,j)$ and yields $x_s \le x_j$ or $x_s > b_j$. Since
$x_s > x_j$, the first alternative fails, so $x_s > b_j$.
\end{proof}

\begin{lemma}\label[lemma]{telescoping}
Fix $i \in [n]$, and let $T$ be any nonempty subset of
$\{\, r : r \pgt{\vx} i \,\}$, with largest element $r_p := \max (T)$. For $\vx$ with conditions \ref{cond:C1} and \ref{cond:C2},
\[
  x_i + \sum_{r \in T} (x_r - b_r) > x_{r_p}.
\]
\end{lemma}

\begin{proof}
Write $T = \{r_1 < r_2 < \cdots < r_p\}$ and set $r_0 := i$. For each
$1 \le l \le p$ we have $r_{l-1} \le_{\vx} r_l$: by $r_l \in T$ we have $r_l \pgt{\vx} i$, and since $i \le r_{l-1} < r_l$, left closure of the position order gives $r_{l-1} \le_{\vx} r_l$. \cref{claim} then yields
$x_{r_{l-1}} > b_{r_l}$. Grouping the sum accordingly,
\[
  x_i + \sum_{l=1}^{p} (x_{r_l} - b_{r_l})
    = \sum_{l=1}^{p} \big(x_{r_{l-1}} - b_{r_l}\big) + x_{r_p}
    > x_{r_p},
\]
since every bracketed term is positive.
\end{proof}

We now formally show the combinatorial characterization of ID vectors. Note that $\vk_{\va, \vb}$ and $k_i$ is used whenever we deal with a sequence that by assumption is a valid ID vector for a given $\vb$ and $\va \in \lowid{\vb}$. Conversely, we use $\vx$ and $x_i$ to denote a sequence with the proposed characterization but which we have yet to prove to be a valid ID vector under $\vb$.

\begin{proposition}
\label[proposition]{k-comb-char}
A vector $\vx$ is equal to $\vk_{\va,\vb}$ for a unique $\va \in \lowid{\vb}$ if and only if $\vx$ satisfies conditions \ref{cond:C1} and \ref{cond:C2}.

\end{proposition}

\begin{proof}

\medskip
Suppose $\va \in \lowid{\vb}$ and let $\vx = \vk_{\va,\vb}$. For condition~\ref{cond:C1}, note
that $\vm_{\va,\vb} \ge 0$; moreover, since
$b_i - a_i = \sum_{r \pge{\va} i} m_r \ge m_i$ and $\va \ge 0$, we have
$\vm_{\va,\vb} \le \vb - \va \le \vb$. Hence
$0 \le \vk_{\va,\vb} = \vb - \vm_{\va,\vb} \le \vb$.

For condition~\ref{cond:C2}, fix $i < j$ with $b_r \ge b_j$ for all $i < r < j$, and assume
$k_i \le b_j$, since otherwise the second alternative already holds. Beginning with $j = j_0$, we construct indices $j_0 > j_1 > j_2 > \cdots \ge i$ with $k_{j_0} \ge k_{j_1} \ge k_{j_2} \ge \cdots \ge k_i$.
Suppose $j_s$ has been defined with $i \le j_s \le j$. If $j_s = i$, we stop.
Otherwise $i < j_s$ and $k_i \le b_{j_s}$: for $s = 0$ this is our assumption, and for
$s > 0$ we have $i < j_s < j$, so $b_{j_s} \ge b_j \ge k_i$. By the contrapositive of \cref{lem:k-greater-b}, $i \nleq_{\va} j_s$. Let $j_{s+1}$ be the largest index
$q < j_s$ with $q \nleq_{\va} j_s$; it exists and satisfies $j_{s+1} \ge i$, since $i$
is such an index. By \cref{lem:left-boundary}, $k_{j_{s+1}} \le k_{j_s}$. The indices
strictly decrease and never pass below $i$, so the process stops at some $j_m = i$, giving $k_i = k_{j_m} \le \cdots \le k_{j_1} \le k_j.$ Thus \ref{cond:C2} holds.

\medskip
Conversely, suppose $\vx$ satisfies \ref{cond:C1} and \ref{cond:C2}. We reconstruct a candidate $\va$ from $\vx$, show that their position orders agree, and verify that $\va\le_T \vb$. We define the reconstruction formula as

$$a_i = x_i + \sum_{r \pgt{\vx} i} (x_r - b_r).$$
Note that, considering $\vk_{\va, \vb}$ in place of $\vx$, this is equivalent to \eqref{k_breakdown} if we show that $i \le_{\va} j \Longleftrightarrow i \le_{\vx} j$. 

Let $\vx$ be a vector with conditions \ref{cond:C1} and \ref{cond:C2}. Assume $i \le_{\vx} j$, and assume also $i<j$ since $i \le_{\va} i \Longrightarrow i \le_{\vx} i$ is trivial. Necessarily, for $r > j$ we have $r \pgt{\vx} i \Longleftrightarrow r \pgt{\vx} j$, by transitivity and left-closure of $\le_{\vx}$. We employ the telescoping bound from \cref{telescoping} to simplify the sum of $x_r -b_r$ over those $r \pgt{\vx} i$ with $r \le j$. There will always be at least one such $r$, specifically $r_p = j$. Together, we have

\begin{align*}
a_i &= x_i + \sum_{r \pgt{\vx} i} (x_r - b_r) \\
    &= x_i + \sum_{\substack{r \pgt{\vx} i \\ r \le j}} (x_r - b_r) + \sum_{r \pgt{\vx} j} (x_r - b_r) \\
    &> x_j + \sum_{r \pgt{\vx} j} (x_r - b_r) \\
    &= a_j.
\end{align*}
Hence, $i \le_{\vx} j \Longrightarrow a_i > a_j$. Moreover, since $i \le_{\vx} j$ implies $s\le_{\vx} j$ for every $i \le s < j$ due to left-closure, it also follows that $a_{i+1}, \ldots, a_{j-1} > a_j$. Thus, $i\le_{\vx}j \Longrightarrow i \le_{\va} j$.

We show by similar means $i \nleq_{\vx} j \Longrightarrow i \nleq_{\va} j$, whose contrapositive gives the desired $i \le_{\va} j \Longrightarrow i \le_{\vx} j$. Assume $i \nleq_{\vx} j$ and choose maximal $t$ with $x_t \le x_j$, $i \le t < j$. If there exist no $r$ with $r \pgt{\vx} j$ and $x_r \ge x_t$, then the assumption $x_j \ge x_t$ immediately gives the inequality below. Otherwise, the set of viable $r$ is nonempty and we use \cref{telescoping} to show

$$x_j + \sum_{\substack{r \pgt{\vx} j \\ x_r \ge x_t}} (x_r - b_r) > x_{r_p} \ge x_t,$$
where the last step follows from the summation condition $x_r \ge x_t$. 

Given the conditions of $t$ and $j$ established thus far, for $r > t$, we assert that \; $r \pgt{\vx} t \Longleftrightarrow (r \pgt{\vx} j$ and $x_r < x_t)$. First assume $r \pgt{\vx} t$. This immediately gives the second condition on the right-hand side. Note for the first condition that no $r \le j$ satisfies $r \pgt{\vx} t$, due to $x_{t+1}, \ldots, x_{j-1} > x_j \ge x_t$, so $t < j < r$ necessarily and left-closure within the partial ordering gives $r \pgt{\vx} j$. Now assume $r \pgt{\vx} j$ and $x_r < x_t$. By the maximality of $t$ we know $x_{t+1}, \ldots, x_{j-1} > x_j > x_r$, and paired with the given $x_t > x_r$ we reach $r \pgt{\vx} t$. We use this fact to rewrite the $r \pgt{\vx} t$ summation condition in terms of $j$.

\begin{align*}
a_t &= x_t + \sum_{r \pgt{\vx} t} (x_r - b_r) \\
    &= x_t + \sum_{\substack{r \pgt{\vx} j \\ x_r < x_t}} (x_r - b_r) \\
    &\le x_j + \sum_{\substack{r \pgt{\vx} j \\ x_r \ge x_t}} (x_r - b_r) + \sum_{\substack{r \pgt{\vx} j \\ x_r < x_t}} (x_r - b_r) \\
    &= a_j.
\end{align*}
The existence of some $a_t \le a_j$ with $i \le t < j$ means $i \nleq_{\va} j$. Hence, $i \nleq_{\vx} j \Longrightarrow i \nleq_{\va} j$, and altogether we have $i \le_{\vx} j \Longleftrightarrow i \le_{\va} j$; from the initial reconstruction formula and this result, we have

$$a_i = x_i + \sum_{r \pgt{\vx} i} (x_r - b_r) = x_i + \sum_{r \pgt{\va} i} (x_r - b_r),$$
returning the original definition \eqref{k_breakdown} of $\vk_{\va, \vb}$. 

However, we must still prove that the reconstructed $\va$ is in $\lowid{\vb}$ to prove $\vx = \vk_{\va, \vb}$, since $\vk_{\va, \vb}$ is defined on a given $\va \in \lowid{\vb}$. $\va \le \vb$ follows from \ref{cond:C1}, since $\va$ is equivalent to $\vx \le \vb$ with non-positive $x_r-b_r$ terms added to it. If there are no $r$ with $r \pgt{\vx} i$, then $a_i = x_i \ge 0$. Otherwise, we use the telescoping bound from \cref{telescoping} to show
$$a_i = x_i + \sum_{r \pgt{\vx} i} (x_r - b_r) > x_{r_p} \ge 0.$$
From here, we show that whenever $i < j$ and $i \le_{\va} j$, the stable Tamari condition holds. Note that the $i=j$ case is trivial. Assuming $i \le_{\va} j$, we also have from the earlier result $i \le_{\vx} j \Longleftrightarrow i \le_{\va} j$ that $i \le_{\vx} j$. Then, using our reconstruction formula for $\va$,

\begin{align*}
a_i - a_j &= x_i + \sum_{r \pgt{\vx} i} (x_r - b_r) - \Big(x_j + \sum_{r \pgt{\vx} j} (x_r - b_r)\Big) \\
          &= x_i - x_j + \sum_{\substack{r \pgt{\vx} i \\ r \le j}} (x_r - b_r) \\
          &= b_i - b_j + \sum_{\substack{r \pge{\vx} i \\ r < j}} (x_r - b_r) \\
          &\le b_i - b_j.
\end{align*}

Thus, $\va \in \lowid{\vb}$. Now set $\vm=\vb-\vx$; having proven that the position orders of $\va$ and $\vx$ agree, the reconstruction formula gives $\vb -\va=\vm M_{\mathbf a}$. Therefore, $\vm=(\vb-\va)M_{\va}^{-1}=\vm_{\va,\vb}$, and hence $\vx=\vb-\vm=\vb-\vm_{\va, \vb}=\vk_{\va,\vb}$. Finally, $\va$ is unique: if $\va'' \in \lowid{\vb}$ also satisfies $\vk_{\va'', \vb} = \vx$, then $\vm_{\va'',\vb} = \vb - \vx = \vm_{\va,\vb}$, so $\va'' = \va$ by the injectivity of $\va \mapsto \vm_{\va,\vb}$ on $\lowid{\vb}$. Thus, for any vector $\vx$ with conditions \ref{cond:C1} and \ref{cond:C2}, $\vx = \vk_{\va, \vb}$ for a unique $\va \in \lowid{\vb}$.

\end{proof}

\section{Level components and the bijective map}\label{sec-3}

\subsection{The level framework}
Condition \ref{cond:C2} of the characterization of ID vectors depends on the entries of $\vb$ between two positions. We therefore organize these positions into level sets, whose interval components provide a convenient way to describe the reversal bijection. For $u \in \N$, define \emph{level u} of $\vb$ as
\[
  L_u(\vb) := \{\, i \in [n] : b_i \ge u \,\}.
\]
This is a disjoint union of intervals. We call each interval component of
$L_u(\vb)$ a \emph{level-$u$ component} and denote it $P$.

\begin{example}
For the running example $\vb = (3,5,2,4,1,4,3)$, the level sets decompose as
\begin{align*}
  L_5(\vb) &= [2,2], \\
  L_4(\vb) &= [2,2] \cup [4,4] \cup [6,6], \\
  L_3(\vb) &= [1,2] \cup [4,4] \cup [6,7], \\
  L_2(\vb) &= [1,4] \cup [6,7], \\
  L_1(\vb) &= [1,7].
\end{align*}
The entries strictly smaller than $u$ separate the components; at level $3$, for instance,
$b_5 = 1$ splits $[4,4]$ from $[6,7]$, and $b_3 = 2$ splits $[1,2]$ from $[4,4]$.
\end{example}

If $P = [p, q]$ is a level-$u$ component of $\vb$, let
  $P^* := [\,n+1-q,\ n+1-p\,]$
be the corresponding level-$u$ component of $\rev(\vb)$. For each level $u$ and component $P \subseteq L_u(\vb)$, set
$
  c_u(P) := \big|\{\, i \in P : k_i = u \,\}\big|,
$
meaning the number of copies of $u$ in component $P$. These definitions will help to express how the bijection rearranges the entries in the ID vector. However, we first provide a more convenient form of the characterization of ID vectors using the level framework.

\begin{proposition}
\label[proposition]{new-k-char}
Assuming condition \ref{cond:C1}, condition \ref{cond:C2} of ID vectors is equivalent to the following: If both $i$ and $j$ with $i < j$ lie in the same level-$u$ component and $k_i=u$, then it must be that $k_j \geq u$.
\end{proposition}

\begin{proof}

Assume first that the new characterization holds. Suppose $i<j$ is a pair considered in \ref{cond:C2}, meaning $b_r \ge b_j$ for $i < r < j$, and suppose $k_i \leq b_j$. Set $u = k_i$. Both $i$ and $j$ must lie in level $u$ since $b_i$, $b_j \ge k_i = u$. Additionally, since $b_j \geq u$ and all intermediate $b_r \geq b_j \geq u$, the entire interval from i to j lies in the same level-$u$ component. The local condition from the new characterization then gives $k_j \geq u = k_i$, which is exactly the first alternative in condition \ref{cond:C2} of the original characterization.

Next, assume that the original characterization holds. Suppose however that the new characterization does not, and choose $i < q$ in $P$ with
$k_i = u$, $k_q < u$, and $q - i$ minimal. Since $q \in P$ we have
$b_q \ge u = k_i$, so $k_i \le b_q$, while $k_i = u > k_q$. If the pair $(i,q)$
is considered --- that is, $b_r \ge b_q$ for all $i < r < q$ --- then $k_i > k_q$
and $k_i \le b_q$ violate condition \ref{cond:C2}.
Otherwise let $t$ be the largest index with $i < t < q$ and $b_t < b_q$. By
maximality, $b_r \ge b_q$ for all $t < r < q$, so $(t,q)$ is considered; and
$k_t \le b_t < b_q$, so condition \ref{cond:C2} forces $k_t \le k_q < u$. As $P$ is an
interval and $i < t < q$, we have $t \in P$, so $(i,t)$ is a pair in $P$ with
$k_i = u > k_t$ and $t - i < q - i$, contradicting minimality of our original choice of $(i,q)$.

\end{proof}

\subsection{The bijection}

We now define the bijective map. Given an $\va \in \lowid{\vb}$, the bijection takes the associated ID vector and constructs a corresponding ID vector for some $\va' \in \lowid{\rev(\vb)} \;$, utilizing the positional information given by level sets in the process.

\begin{definition}

We define $\Phi_{\vb}(\vk_{\va, \vb}) = \vy$ by processing the integer levels
$u = \max(k_i), \max(k_i)-1, \dots, 1, 0$. For each component $P \subseteq
L_u(\vb)$, place $c_u(P)$ copies of $u$ into the rightmost currently empty
positions of $P^*$. Components at a fixed level are disjoint, so they may be processed simultaneously.

\end{definition}

\begin{example}
We compute $\vy = \Phi_{\vb}(\vk_{\va, \vb})$ for the running example, with
$\vk_{\va, \vb} = (0,5,1,3,1,4,3)$ and $\rev(\vb) = (3,4,1,4,2,5,3)$. Processing
$u = 5,4,\dots,0$ and listing only the components with $c_u(P) > 0$:
\[
\renewcommand{\arraystretch}{1.25}
\begin{array}{c|c|c|c|c}
  u & P & c_u(P) & P^* & \text{$\vy$ after placing} \\ \hline
  5 & [2,2] & 1 & [6,6] & (\cdot,\cdot,\cdot,\cdot,\cdot,5,\cdot) \\
  4 & [6,6] & 1 & [2,2] & (\cdot,4,\cdot,\cdot,\cdot,5,\cdot) \\
  3 & [4,4] & 1 & [4,4] & (\cdot,4,\cdot,3,\cdot,5,\cdot) \\
  3 & [6,7] & 1 & [1,2] & (3,4,\cdot,3,\cdot,5,\cdot) \\
  1 & [1,7] & 2 & [1,7] & (3,4,\cdot,3,1,5,1) \\
  0 & [1,7] & 1 & [1,7] & (3,4,0,3,1,5,1)
\end{array}
\]
The result is $\vy = \Phi_{\vb}(\vk_{\va, \vb}) = (3,4,0,3,1,5,1)$. We validate that the conditions for ID vectors under $\rev(\vb)$ hold here. Specifically, it is $\vk_{\va',\rev(\vb)}$ for $\va' = (2,3,0,2,1,3,1) \in \lowid{\rev(\vb)}$, which we find through
$\rev(\vb) - \vk_{\va', \rev(\vb)} = \vm_{\va',\rev(\vb)} = (0,0,1,1,1,0,2)$ and \eqref{m}. Since $|\vb| = |\rev(\vb)| = 22$ and $|\vk_{\va, \vb}| = |\vk_{\va', \rev(\vb)}| = 17$, $\hgt(\va',\rev(\vb)) = 22 - 17 = 5 = \hgt(\va,\vb)$.
\end{example}

In the stated order, we now prove that the map is well-defined, constructs a valid ID vector under $\rev(\vb)$, and is an involution.

\begin{proposition}
For every ID vector $\vk_{\va, \vb}$ under $\vb$, the construction of $\Phi_{\vb}(\vk_{\va, \vb})=\vy$ can be completed.
\end{proposition}

\begin{proof}
Fix a level $u$ and a component $P \subseteq L_u(\vb)$. Before level $u$ is
processed, only values greater than $u$ have been placed.

If $i \in P$ and $k_i = v > u$, then the level-$v$ component containing $i$ is
contained in $P$, because $L_v(\vb) \subseteq L_u(\vb)$. Hence the image of that entry was placed in $P^*$. Furthermore, if an entry $k_i = v > u$ is placed in $P^*$, then we know that it must come from a level-$v$ component contained within $P$, again since $L_v(\vb) \subseteq L_u(\vb)$. All together, this means that any $v$ within $P$ is placed in $P^*$ and no $v$ from outside $P$ can be placed in $P^*$. Thus, the number of already occupied positions of $P^*$ is exactly $\big|\{\, i \in P : k_i > u \,\}\big|.$
The number of empty positions is therefore
\[
  |P| - \big|\{\, i \in P : k_i > u \,\}\big|
  = \big|\{\, i \in P : k_i \le u \,\}\big|
  \ge \big|\{\, i \in P : k_i = u \,\}\big|
  = c_u(P).
\]
Hence all copies of $u$ can be placed.
\end{proof}

We now show that the map constructs a valid ID vector under $\rev(\vb)$. To do so, we show that the conditions that define ID vectors under $\rev(\vb)$ apply to $\vy = \Phi_{\vb}(\vk_{\va, \vb})$.

\begin{proposition}
\label[proposition]{valid-bij}
If $\vk_{\va, \vb}$ is a valid ID vector under $\vb$, then $\vy = \Phi_{\vb}(\vk_{\va, \vb})$ is a valid ID vector under $\rev(\vb)$.
\end{proposition}

\begin{proof}
We show that both conditions \ref{cond:C1} and \ref{cond:C2} hold for the constructed $\vy$. Condition \ref{cond:C1} is immediate: $\vy$ is a rearrangement of the terms of $\vk_{\va, \vb} \ge 0$, where each $y_i = u$ has been placed in a level-$u$ component with $\rev(\vb)_i \ge u = y_i$, so $0 \le \vy \le \rev(\vb)$.

We show condition \ref{cond:C2} in its new form from \cref{new-k-char}. Fix a level-$u$ component $P$. The only entries that have been placed so far in $P^*$ are greater than $u$. All $c_u(P)$ copies of $u$ in $P$ are now placed in the rightmost available positions within $P^*$. This means that for $i,j \in P^*$, $i<j$, $y_i = u$ necessarily implies $y_j \ge u$. In other words, if any entry $y_j$ in $P^*$ is to the right of another entry $y_i=u$ in $P^*$, it must have been placed either before or simultaneously, hence $y_j \ge u$. By \cref{new-k-char}, we now also have condition \ref{cond:C2} of the original ID vector characterization given by \cref{k-comb-char}. Thus, the resulting sequence $\vy$ is $\vk_{\va', \rev(\vb)}$ for some $\va' \in \lowid{\rev(\vb)}$, or equivalently is a valid ID vector under $\rev(\vb)$.

\end{proof}

\begin{proposition}
\label[proposition]{involution}
The maps $\Phi_{\vb}$ and $\Phi_{\rev(\vb)}$ are inverse to one another.
\end{proposition}

\begin{proof}
It suffices to show $\Phi_{\rev(\vb)}\big(\Phi_{\vb}(\vk_{\va, \vb})\big) = \vk_{\va, \vb}$ for every ID vector $\vk_{\va, \vb}$ of $\vb$; interchanging $\vb$ and $\rev(\vb)$ then gives the other
composition. Write $\vk_{\va', \rev(\vb)} = \Phi_{\vb}(\vk_{\va, \vb})$.

\medskip
Apply $\Phi_{\rev(\vb)}$ to $\vk_{\va', \rev(\vb)}$, processing the levels
$u = \max_i (k_i), \dots, 1, 0$ in descending order, and assume inductively that
every value greater than $u$ has already been returned to its original
position. Fix a level-$u$ component $P$ of $\vb$, with mirror $P^*$ a level-$u$
component of $\rev(\vb)$. In the forward map, exactly $c_u(P)$ copies of $u$ were
placed into $P^*$; these are precisely the entries of $\vk_{\va', \rev(\vb)}$ equal to $u$ in
$P^*$, so $\Phi_{\rev(\vb)}$ now returns exactly $c_u(P)$ copies of $u$ to $P$. By
the induction hypothesis the positions of $P$ holding values greater than $u$
are already refilled in their original spots, so the empty positions of $P$ are
exactly those whose original value is at most $u$. The map $\Phi_{\rev(\vb)}$
places the $c_u(P)$ copies of $u$ into the rightmost of these empty positions; since $\vk_{\va,\vb}$ is an ID vector under $\vb$, \cref{new-k-char} shows that these are precisely the positions at which $\vk_{\va,\vb}$ has its copies of $u$. Hence, every level-$u$ entry returns to its original position, and descending on $u$ gives
\[
  \Phi_{\rev(\vb)}\big(\Phi_{\vb}(\vk_{\va, \vb})\big) = \vk_{\va, \vb}.
\]
Interchanging $\vb$ and $\rev(\vb)$ gives the other composition.
\end{proof}

We now use the established bijection to prove the reversal symmetry of $f_{\vb}(t)$.

\begin{theorem}\label[theorem]{main-result}
For every $\vb \in \N^n$,
\[
  f_{\vb}(t) = f_{\rev(\vb)}(t).
\]
\end{theorem}

\begin{proof}
The map $\Phi_{\vb}$ uses sequences of the form $\vk_{\va, \vb}$ to construct sequences of the form $\vk_{\va', \rev(\vb)}$. Since it only rearranges the entries of $\vk_{\va, \vb}$,
\[
  |\vk_{\va, \vb}| = |\vk_{\va', \rev(\vb)}|.
\]
By \cref{valid-bij,involution}, $\Phi_{\vb}$ induces a bijection
\[
  \lowid{\vb} \to \lowid{\rev(\vb)}, \qquad \va \mapsto \va'.
\]
Since $|\vb| = |\rev(\vb)|$, \eqref{k_height} gives
\[
  \hgt(\va', \rev(\vb)) = |\rev(\vb)| - |\vk_{\va', \rev(\vb)}| = |\vb| - |\vk_{\va, \vb}| = \hgt(\va,\vb).
\]
Thus, the induced bijection preserves height, and summing $t^{\hgt}$ over the two lower ideals proves the result.
\end{proof}

\section*{Acknowledgments}
I thank my research advisor, Anna Pun, for her generous feedback and support throughout this project; though I initially developed ID vectors with a visual conception of their characterization, she helped me to formalize them in \eqref{k_breakdown} via the height matrix from her paper \cite{chau2026stabletamarilattice}. I am also grateful to the high school research program of the Discrete Mathematics Group at Baruch College, which gave me the foundation to begin working on open problems.

Claude (Anthropic) was used as an interactive aid during this project, including assistance with figure generation, computations for verifying conjectures, providing the paper's running example (given my criterion for an ideal example), LaTeX formatting, and general editing feedback. It also gave the \cref{telescoping} telescoping bound idea and suggested the $q-i$ minimality condition in the \cref{new-k-char} proof. However, I constructed the entire argument for the proofs of \cref{k-comb-char,new-k-char}, using these tools as necessary along the way. Except for the assistance described above and results cited from other sources, the arguments are my own.

\bibliographystyle{amsplain}
\bibliography{refs}

\providecommand{\bysame}{\leavevmode\hbox to3em{\hrulefill}\thinspace}
\providecommand{\MR}{\relax\ifhmode\unskip\space\fi MR }
\providecommand{\MRhref}[2]{%
  \href{http://www.ams.org/mathscinet-getitem?mr=#1}{#2}
}
\providecommand{\href}[2]{#2}
\begin{thebibliography}{1}

\bibitem{MR3029440}
Fran{\c{c}}ois Bergeron and Louis-Fran{\c{c}}ois Pr\'eville-Ratelle, \emph{Higher trivariate diagonal harmonics via generalized {T}amari posets}, J. Comb. \textbf{3} (2012), no.~3, 317--341. \MR{3029440}

\bibitem{MR2880681}
Mireille Bousquet-M\'elou, \'Eric Fusy, and Louis-Fran{\c{c}}ois Pr\'eville-Ratelle, \emph{The number of intervals in the {$m$}-{T}amari lattices}, Electron. J. Combin. \textbf{18} (2011), no.~2, Paper 31, 26. \MR{2880681}

\bibitem{chau2026stabletamarilattice}
Herman Chau, Spencer Daugherty, Sean Grate, Mark~D. Haiman, Juan Carlos~Mart{\'\i}nez Mori, and Anna Pun, \emph{The stable {T}amari lattice}, 2026, Preprint, arXiv:2610.39257.

\bibitem{MR3632566}
Louis-Fran{\c{c}}ois Pr\'eville-Ratelle and Xavier Viennot, \emph{The enumeration of generalized {T}amari intervals}, Trans. Amer. Math. Soc. \textbf{369} (2017), no.~7, 5219--5239. \MR{3632566}

\end{thebibliography}

\end{document}